\documentclass[12pt]{amsart}
\usepackage{amsfonts}
\usepackage{amssymb,mathrsfs}
\usepackage{amsmath}
\usepackage{enumerate}
\usepackage{latexsym}
\usepackage{graphicx}
\usepackage{bm}
\usepackage{subfigure}
\usepackage{float}
\usepackage{amsmath}
\usepackage{amsthm}
\usepackage{verbatim}
\usepackage{vmargin}
\usepackage{amstext}
\usepackage{array}
\usepackage{booktabs}
\usepackage{indentfirst}
\usepackage[compress]{cite}
\usepackage{url}
\usepackage{geometry}
\numberwithin{equation}{section}

\newtheorem{thm}{Theorem}[section]

\newtheorem{lma}[thm]{Lemma}

\newtheorem{defn}[thm]{Definition}
\newtheorem{rem}[thm]{Remark}

\theoremstyle{definition}

\begin{document}
\setlength{\oddsidemargin}{80pt}
\setlength{\evensidemargin}{80pt}
\setlength{\topmargin}{80pt}
\author{Jun~Xian, Xiaoda~Xu}

\address{J.~Xian\\Department of Mathematics and Guangdong Province Key Laboratory of Computational Science
 \\Sun
Yat-sen University
 \\
510275 Guangzhou\\
China.} \email{xianjun@mail.sysu.edu.cn}

\address{X.~Xu\\Department of Mathematics
 \\Sun
Yat-sen University
 \\
510275 Guangzhou\\
China.} \email{xuxd26@mail2.sysu.edu.cn}

\title[Expected star discrepancy for stratified sampling]{On the lower expected star discrepancy for jittered sampling than simple random sampling}

\keywords{Expected star discrepancy; Stratified sampling.}

\date{\today}

\pdfoutput=1

\subjclass[2010]{ 65C10,  11K38,  65D30. }
\begin{abstract}
We compare expected star discrepancy under jittered sampling with simple random sampling, and the strong partition principle for the star discrepancy is proved.
\end{abstract}

\maketitle
\section{Introduction}\label{intro}

Classical jittered sampling (JS) acquires better expected star discrepancy bounds than using traditional Monte Carlo (MC) sampling, see \cite{Doerr2}. This means jittered sampling is the refinement of the traditional Monte Carlo method, the problem now is that is there a direct comparison of star discrepancy itself, not a bound comparison. This involves the strong partition principle for the star discrepancy version, see an introduction in \cite{KP2}. Among various techniques for measuring the irregularities of point distribution, star discrepancy is the most efficient, which is well established and has found applications in various areas such as computer graphics, machining learning, numerical integration, and financial engineering, see \cite{Owen, MLR, Dick2013QMC, Pgm}.

\textbf{Star discrepancy}. The star discrepancy of a sampling set $P_{N,d}=\{t_{1},t_{2},\ldots,t_{N}\}$ is defined by:

\begin{equation}\label{*-d}D_{N}^{*}\left(t_{1}, t_{2}, \ldots, t_{N}\right):=\sup _{x \in[0,1]^{d}}\left|\lambda([0, x])-\frac{\sum_{n=1}^{N} I_{[0, x]}\left(t_{n}\right)}{N}\right|,\end{equation} where $\lambda$ denotes the $d$-dimensional Lebesgue measure and $I_{[0,x]}$ denotes the characteristic function defined on the $d$-dimensional rectangle $[0,x]$.

The special deterministic point set constructions are called \textbf{low discrepancy point sets}, the best known asymptotic upper bounds for star discrepancy are of the form $$O(\frac{(\ln N)^{\alpha_{d}}}{N}),$$ where $\alpha_{d}\ge 0$ for fixed $d$. Studies on such point sets, see \cite{HNie,Dick2010QMC}. We introduce random factors in the present paper. Formers have conducted extensive research on the comparisons between expected discrepancy. In \cite{KP2}, it is proved that jittered sampling constitutes a set whose expected $L_p-$discrepancy is smaller than that of purely random points. Further, a theoretical conclusion that the jittered sampling does not have the minimal expected $L_2-$discrepancy is presented in \cite{KP}. We study expected star discrepancy under different partition models, which are the following, 

\begin{equation}\label{formula12}
\mathbb{E}(D_{N}^{*}(Y))< \mathbb{E}(D_{N}^{*}(X)),
\end{equation}
where $X$ denotes a simple random sampling set, $Y$ denotes a stratified sampling point set that is uniformly distributed in the grid-based stratified subsets. Since 2016 in \cite{jittsamp}, F. Pausinger and S. Steinerberger have given the upper and lower bounds for expected star discrepancy under jittered sampling, and they proved the strong partition principle for $L_2-$discrepancy, they mentioned whether this conclusion could be generalized to star discrepancy. Then, in 2021 \cite{KP}, M. Kiderlen and F. Pausinger referred again to the proof of the strong partition principle for the star discrepancy.

The rest of this paper is organized as follows. Section 2 presents preliminaries on stratified sampling and $\delta-$covers. Section \ref{sec3} presents our main results, which provide comparisons of the expected star discrepancy for simple random sampling and stratified sampling under grid-based equivolume partition. Section \ref{pfmr} includes the proofs of the main result. Finally, in section \ref{conclu} we conclude the paper with a short summary.

\section{Preliminaries on stratified sampling and $\delta-$covers}\label{prelim}

Before introducing the main result, we list the preliminaries used in this paper.

\subsection{Jittered sampling}
\begin{figure*}[h]
\centering
\subfigure[jittered sampling in two dimension]{
\begin{minipage}{7cm}
\centering
\includegraphics[width=0.6\textwidth]{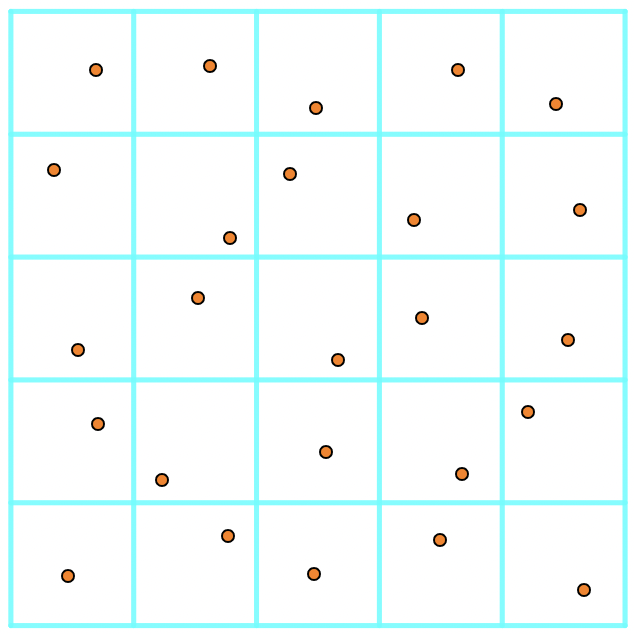}
\end{minipage}
}
\subfigure[jittered sampling in three dimension]{
\begin{minipage}{7cm}
\centering
\includegraphics[width=0.7\textwidth]{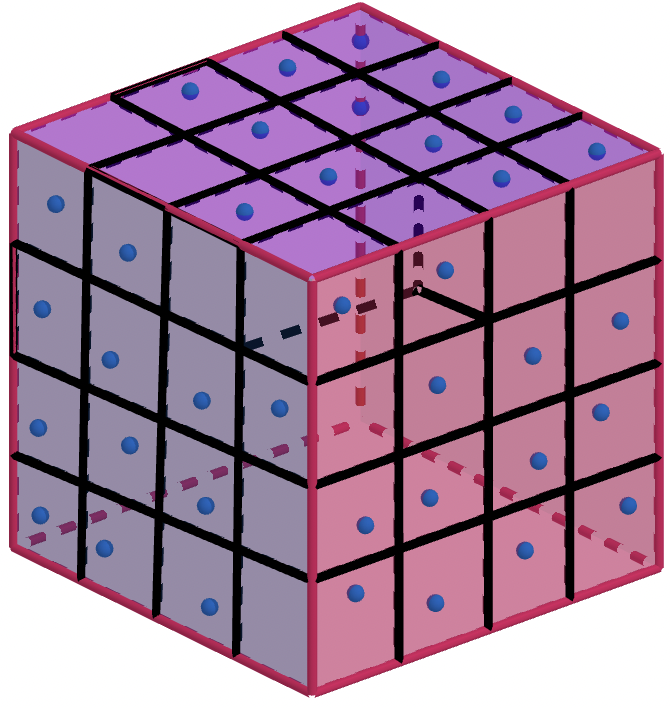}
\end{minipage}
}
\caption{\label{ss0} jittered sampling formed by isometric grid partition.}
\end{figure*}

Jittered sampling is a type of grid-based equivolume partition. $[0,1]^{d}$ is divided into $m^d$ axis parallel boxes $Q_{i},1\leq i\leq N,$ each with sides $\frac{1}{m},$ see illustration of Figure \ref{ss0}. Research on the jittered sampling are extensive, see \cite{shirely1994,Doerr2,KP,KP2,jittsamp}.

For jittered sampling, we consider a rectangle $R=[0,x)$ (we shall call it the test set in the following) in $[0,1]^{d}$ anchored at $0$. For the corresponding isometric grid partition $\Omega=\{Q_{1},Q_{2},
\ldots,Q_{N}\}$ of $[0,1]^{d}$, if we put

$$
I_{N}:=\{j:\partial R \cap Q_{j}\neq \emptyset\},
$$
and

$$
C_{N}:=|I_{N}|,
$$
which means the cardinality of the index set $I_{N}$. Then for $C_N$, easy to show that

\begin{equation}\label{CNbd0}
C_{N}\leq d\cdot N^{1-\frac{1}{d}}.
\end{equation}

\subsection{$\delta-$covers}
Secondly, to discretize the star discrepancy, we use the definition of $\delta-$covers as in \cite{Doerr}, which is well known in empirical process theory, see, e.g., \cite{vawell}.

\begin{defn}
For any $\delta\in (0,1]$, a finite set $\Gamma$ of points in $[0,1)^{d}$ is called a $\delta$-cover of $[0,1)^{d}$, if for every $y\in [0,1)^{d}$, there exist $x, z\in \Gamma\cup \{0\}$ such that $x\leq y\leq z$ and $\lambda([0, z])-\lambda([0, x])\leq \delta$. The number $\mathcal{N}(d,\delta)$ denotes the smallest cardinality of a $\delta$-cover of $[0,1)^{d}$.
\end{defn}

From \cite{Doerr,Gnewuch2008Bracketing}, combining with Stirling's formula, the following estimation for $\mathcal{N}(d,\delta)$ holds, that is, for any $d\ge 1$ and $\delta\in (0,1]$,

\begin{equation}\label{delcovbd}
    \mathcal{N}(d,\delta)\leq 2^d\cdot \frac{e^d}{\sqrt{2\pi d}}\cdot(\delta^{-1}+1)^{d}.
\end{equation}

Let $P=\{p_{1},p_{2},\ldots,p_{N}\}\subset [0,1]^{d}$ and $\Gamma$ be $\delta-$covers, then

$$
D_{N}^{*}(P)\leq D_{\Gamma}(P)+\delta,
$$
where

\begin{equation}\label{dgammap}
D_{\Gamma}(P):=\max_{x\in\Gamma}|\lambda([0,x])-\frac{\sum_{n=1}^{N} I_{[0, x]}\left(p_{n}\right)}{N}|.
\end{equation}

Formula \eqref{dgammap} provides convenience for estimating the star discrepancy.

\section{Expected star discrepancy for stratified and simple random sampling}\label{sec3}

In this section, a comparison of expected star discrepancy for jittered sampling and simple random sampling is obtained, which means the strong partition principle for star discrepancy holds.

\subsection{Comparison of expected star discrepancy under jittered sampling and simple random sampling}

\begin{thm}\label{uniformnoise0}
Let $m,d\in \mathbb{N}$ with $m\ge d\ge 2$. Let $N=m^d$, $X=\{X_{1}, X_{2},\ldots, X_{N}\}$ is a simple random sampling set. Stratified random $d-$dimension point set $Y=\{Y_{1}, Y_{2}, Y_{3}, \ldots, Y_{N}\}$ is uniformly distributed in the grid-based stratified subsets of $\Omega^{*}_{|}$, then

\begin{equation}\label{formula31}
\mathbb{E}(D_{N}^{*}(Y))< \mathbb{E}(D_{N}^{*}(X)).
\end{equation}

\end{thm}

\begin{rem}
The inequality \eqref{formula31} in Theorem \ref{uniformnoise0} presents lower expected star discrepancy under grid-based partition than simple random sampling, which also means [Strong Partition Principle \cite{KP2}] holds, it is stated in \cite{KP} that strong partition principle for the star discrepancy needs to be proved.
\end{rem}

\section{Proofs}\label{pfmr}
In this section, we present the proof of Theorem  \ref{uniformnoise0}. We first introduce Bernstein's inequality, which is a standard result that can be found in textbooks on probability, see, e.g., \cite{FDX2007} and we shall omit the proof here. 

\begin{lma}[\textbf{Bernstein's inequality}]\label{Binequ}
Let $Z_1,\ldots,Z_N$ be independent random variables with expected
values $\mathbb E (Z_j)=\mu_{j}$ and variances $\sigma_{j}^2$ for $j=1,\ldots,N$.  Assume  $|Z_j-\mu_j|\le C$(C is a constant) for each $j$ and set $\Sigma^2:=\sum_{j=1}^N\sigma_{j}^2$, then
for any $\lambda \ge 0$,
$$\mathbb P\left\{\Big|\sum_{j=1}^N [Z_j-\mu_j]\Big|\ge \lambda\right \}\le
2\exp\left(-\frac{\lambda^2}{2\Sigma^2+\frac{2}{3}C \lambda}\right).$$
\end{lma}

\subsection{Proof of Theorem \ref{uniformnoise0}}
For the arbitrary test set $R=[0,x)$, we unify a label $W=\{W_1,W_2,\ldots,W_N\}$ for the sampling point sets formed by different sampling models, and we consider the following discrepancy function,

\begin{equation}\label{dispfunc1}
\Delta_{\mathscr{P}}(x)=\frac{1}{N}
\sum_{n=1}^{N}\mathbf{1}_{R}(W_n)-\lambda(R).
\end{equation}

For the grid-based equivolume partition $\Omega=\{\Omega_1,\Omega_2,\ldots,\Omega_N\}$, we divide the test set $R$ into two parts, one is the disjoint union of $\Omega_{i}$ entirely contained by $R$ and another is the union of remaining pieces which are the intersections of some $\Omega_{j}$ and $R$, i.e.,

\begin{equation}\label{Rtp1}
R=\bigcup_{i\in I_0}\Omega_{i}\cup\bigcup_{j\in J_0}(\Omega_{j}\cap R),
\end{equation}
where $I_0,J_0$ are two index sets.

We set $$T=\bigcup_{j\in J_0}(\Omega_{j}\cap R).$$

Besides, for an equivolume partition $\Omega=\{\Omega_1, \Omega_2, \ldots, \Omega_N\}$ of $[0,1]^{d}$ and the corresponding stratified sampling set $P_{\Omega}$, we have

\begin{equation}\label{for4d4}
\begin{aligned}
\text{Var}(\sum_{n=1}^{N}\mathbf{1}_{R}(P_{\Omega}))&=\sum_{i=1}^{N}\frac{|\Omega_i\cap[0,x]|}{|\Omega_i|}(1-\frac{|\Omega_i\cap[0,x]|}{|\Omega_i|})\\&=N|[0,x]|-N^2\sum_{i=1}^{N}(|\Omega_i\cap[0,x]|)^2.
\end{aligned}
\end{equation}

For sampling sets $Y=\{Y_1,Y_2,\ldots,Y_N\}$ and $X=\{X_1,X_2,\ldots,X_N\}$, we have
\begin{equation}
\text{Var}(\sum_{n=1}^{N}\mathbf{1}_{R}(Y))=N|[0,x]|-N^2\sum_{i=1}^{N}(|\Omega^{*}_{|,i}\cap[0,x]|)^2,
\end{equation}
and 
\begin{equation}
\text{Var}(\sum_{n=1}^{N}\mathbf{1}_{R}(X))=N|[0,x]|-N^2|[0,x]|^2.
\end{equation}

Hence, we have

\begin{equation}\label{varcomp}
\text{Var}(\sum_{n=1}^{N}\mathbf{1}_{R}(Y))\leq \text{Var}(\sum_{n=1}^{N}\mathbf{1}_{R}(X))
\end{equation}

Now, we exclude the equality case in \eqref{varcomp}, which means the following formula holds,

\begin{equation}
N|\Omega^{*}_{|,i}\cap[0,x]|=|[0,x]|,
\end{equation}
for $i=1,2,\ldots,N$ and for almost all $x\in [0,1]^{d}$.

Hence, 

\begin{equation}
\int_{[0,x]}\mathbf{1}_{\Omega^{*}_{|,i}}(y)dy=\int_{[0,x]}\frac{1}{N}dy.
\end{equation}
for almost all $x\in[0,1]^{d}$ and all $i=1,2,\ldots,N$, which implies $\mathbf{1}_{\Omega^{*}_{|,i}}=\frac{1}{N}$, which is not impossible for $N\ge 2.$

For set $T$, we have

\begin{equation}
\text{Var}(\frac{1}{N}
\sum_{n=1}^{N}\mathbf{1}_{T}(Y_n))=\sum_{i=1}^{|J_0|}\frac{|\Omega^{*}_{i,|}\cap[0,x]|}{|\Omega^{*}_{i,|}|}(1-\frac{|\Omega^{*}_{i,|}\cap[0,x]|}{|\Omega^{*}_{i,|}|})
\end{equation}

The same analysis for set $T$ as $R$, we have

\begin{equation}\label{varcom000}
   \text{Var}(\frac{1}{N}
\sum_{n=1}^{N}\mathbf{1}_{T}(Y_n))<\text{Var}(\frac{1}{N}
\sum_{n=1}^{N}\mathbf{1}_{T}(X_n)).
\end{equation}

For test set $R=[0,x)$, we choose $R_0=[0,y)$ and $R_1=[0,z)$ such that $y\leq x\leq z$ and $\lambda(R_1)-\lambda(R_0)\leq \frac{1}{N}$, then ($R_0$,$R_1$) constitute the $\frac{1}{N}-$covers. For $R_0$ and $R_1$, we can divide them into two parts as we did for \eqref{Rtp1} respectively. Let $$T_0=\bigcup_{j\in J_0}(\Omega_{j}\cap R_0),$$ and $$T_1=\bigcup_{j\in J_0}(\Omega_{j}\cap R_1).$$ We have the same conclusion for $T_0$ and $T_1$. In order to unify the two cases $T_0$ and $T_1$ (Because $T_0$ and $T_1$ are generated from two test sets with the same cardinality, and the cardinality is the covering numbers), we consider a set $R'$ which can be divided into two parts

\begin{equation}\label{Rpr1}
R'=\bigcup_{k\in K}\Omega_{k}\cup\bigcup_{l\in L}(\Omega_{l}\cap R'),
\end{equation}
where $K,L$ are two index sets. Moreover, we set the cardinality of $R^{'}\subset[0,1)^{d}$ at most $2^{d-1}\frac{e^{d}}{\sqrt{2\pi d}}(N+1)^{d}$(the $\delta-$covering numbers, where we choose $\delta=\frac{1}{N}$), and we let

$$
T'=\bigcup_{l\in L}(\Omega_{l}\cap R').
$$

We define new random variables $\chi_{j}, 1\leq j\leq |L|$, as follows

$$
    \chi_{j}=\left\{
\begin{aligned}
&1, W_{j}\in \Omega_{j}\cap R',\\&
0, otherwise.
\end{aligned}
\right.
$$
Then,

\begin{equation}\label{murmut0}
\begin{aligned}
&N\cdot D^{*}_{N}\left(W_{1}, W_{2}, \ldots, W_{N};R'\right)\\&=N\cdot D^{*}_{N}\left(W_{1}, W_{2}, \ldots, W_{N};T'\right)\\&=|\sum_{j=1}^{|L|}\chi_{j}-N(\sum_{j=1}^{|L|}\lambda(\Omega_{j}\cap T'))|.
\end{aligned}
\end{equation}

Since

$$
\mathbb{P}(\chi_{j}=1)=\frac{\lambda(\Omega_{j}\cap T')}{\lambda(\Omega_{j})}=N\cdot\lambda(\Omega_{j}\cap T'),
$$
we get

\begin{equation}\label{echi0}
    \mathbf{E}(\chi_{j})=N\cdot\lambda(\Omega_{j}\cap T').
\end{equation}

Thus, from \eqref{murmut0} and \eqref{echi0}, we obtain

\begin{equation}\label{Rprime00}
N\cdot D^{*}_{N}\left(W_{1}, W_{2}, \ldots, W_{N};R'\right)=|\sum_{j=1}^{|L|}(\chi_{j}-\mathbb{E}(\chi_{j}))|.
\end{equation}

Let $$\sigma_{j}^{2}=\mathbb{E}(\chi_{j}-\mathbb{E}(\chi_{j}))^{2}, \Sigma=(\sum_{j=1}^{|L|}\sigma_{j}^{2})^{\frac{1}{2}}.$$

Therefore, from Lemma \ref{Binequ}, for every $R'$, we have,

$$
\mathbb{P}\left
(\Big|\sum_{j=1}^{|L|}(\chi_{j}-\mathbb{E}(\chi_{j}))\Big|>\lambda\right)\leq
2\cdot\exp(-\frac{\lambda^{2}}{2\Sigma^{2}+\frac{2\lambda}{3}}).
$$

Let $\mathscr{B}=\bigcup\limits_{R'}\left(\Big|\sum_{j=1}^{|L|}(\chi_{j}-\mathbb{E}(\chi_{j}))|>\lambda\right),$ then using $\delta-$covering numbers, we have

\begin{equation}\label{pb000}
\mathbb{P}(\mathscr{B})\leq (2e)^{d}\cdot\frac{1}{\sqrt{2\pi d}}\cdot(N+1)^{d}\cdot
\exp(-\frac{\lambda^{2}}{2\Sigma^{2}+\frac{2\lambda}{3}}).
\end{equation}

Combining with \eqref{Rprime00}, we get

\begin{equation}\label{punrp1}
\begin{aligned}
&\mathbb{P}\Big(\bigcup_{R'}\left(N\cdot D_{N}^{*}\left(W_{1}, W_{2}, \ldots, W_{N};R'\right)>\lambda\right)\Big)\\&\leq (2e)^{d}\cdot\frac{1}{\sqrt{2\pi d}}\cdot(N+1)^{d}\cdot
\exp(-\frac{\lambda^{2}}{2\Sigma^{2}+\frac{2\lambda}{3}}).
\end{aligned}
\end{equation}

For point sets $Y$ and $X$, if we let $$\Sigma_{0}^2=\text{Var}(\sum_{n=1}^{N}\mathbf{1}_{T'}(Y_n)), \Sigma_{1}^2=\text{Var}(\sum_{n=1}^{N}\mathbf{1}_{T'}(X_n)).$$

Then \eqref{varcom000} implies

$$\Sigma_{0}^2<\Sigma_{1}^2.$$

Besides, as \eqref{punrp1}, we have

\begin{equation}\label{purpy2}
\begin{aligned}
&\mathbb{P}\Big(\bigcup_{R'}\left(N\cdot D_{N}^{*}\left(Y_{1}, Y_{2}, \ldots, Y_{N};R'\right)>\lambda\right)\Big)\\&\leq (2e)^{d}\cdot\frac{1}{\sqrt{2\pi d}}\cdot(N+1)^{d}\cdot
\exp(-\frac{\lambda^{2}}{2\Sigma_0^{2}+\frac{2\lambda}{3}}),
\end{aligned}
\end{equation}
and

\begin{align*}
&\mathbb{P}\Big(\bigcup_{R'}\left(N\cdot D_{N}^{*}\left(X_{1}, X_{2}, \ldots, X_{N};R'\right)>\lambda\right)\Big)\\&\leq (2e)^{d}\cdot\frac{1}{\sqrt{2\pi d}}\cdot(N+1)^{d}\cdot
\exp(-\frac{\lambda^{2}}{2\Sigma_1^{2}+\frac{2\lambda}{3}}),
\end{align*}
respectively.

Suppose $A(d,q,N)=d\ln(2e)+d\ln(N+1)-\frac{\ln(2\pi d)}{2}-\ln(1-q)$, and we choose

$$
\lambda=\sqrt{2\Sigma_0^2\cdot A(d,q,N)+\frac{A^2(d,q,N)}{9}}+\frac{A(d,q,N)}{3}
$$
in \eqref{purpy2}, then we have

\begin{equation}\label{punionrp}
\mathbb{P}\Big(\bigcup_{R'}\left(N\cdot D_{N}^{*}\left(Y_{1}, Y_{2}, \ldots, Y_{N};R'\right)>\lambda\right)\Big)\leq 1-q.
\end{equation}

Hence, from \eqref{punionrp}, it can easily be verified

\begin{equation}\label{maxdnstar}
\max_{R_i,i=0,1}
D_{N}^{*}(Y_{1},Y_{2},\ldots,Y_{N};R_{i})\leq \frac{\sqrt{2\Sigma_0^2\cdot A(d,q,N)+\frac{A^2(d,q,N)}{9}}}{N}+\frac{A(d,q,N)}{3N}
\end{equation}
holds with probability at least $q$.

From \eqref{dgammap}, combining with $\delta-$covering numbers(where $\delta=\frac{1}{N}$), we get,

\begin{equation}\label{dnsy1}
\begin{aligned}
    D_{N}^{*}(Y)&\leq \frac{\sqrt{2\Sigma_0^2\cdot A(d,q,N)+\frac{A^2(d,q,N)}{9}}}{N}+\frac{A(d,q,N)+3}{3N}\\&
    \leq (\sqrt{2}\cdot\Sigma_0+1)\frac{A(d,q,N)}{N}
\end{aligned}
\end{equation}

holds with probability at least $q$, the last inequality in \eqref{dnsy1} holds because $A(d,q,N)\ge 3$ holds for all $q\in (0,1)$.

Same analysis with point set $X$, we have

\begin{equation}\label{DNsgm1}
\begin{aligned}
    D_{N}^{*}(X)&\leq \frac{\sqrt{2\Sigma_1^2\cdot A(d,q,N)+\frac{A^2(d,q,N)}{9}}}{N}+\frac{A(d,q,N)+3}{3N}\\&
    \leq (\sqrt{2}\cdot\Sigma_1+1)\frac{A(d,q,N)}{N}
\end{aligned}
\end{equation}

holds with probability at least $q$.

Now, we fix a probability value $q_0\in(0,1)$ in \eqref{dnsy1}, i.e., we suppose \eqref{dnsy1} holds with probability exactly $q_0$, where $q_0\in [q,1)$. Choose this $q_0$ in \eqref{DNsgm1}, we have

$$
    D_{N}^{*}(X)\leq (\sqrt{2}\cdot\Sigma_1+1)\frac{A(d,q_0,N)}{N},
$$
holds with probability $q_0.$

Therefore from $\Sigma_0<\Sigma_1$, we obtain,

\begin{equation}\label{dnsxq1}
     D_{N}^{*}(X)\leq (\sqrt{2}\cdot\Sigma_0+1)\frac{A(d,q_0,N)}{N}
\end{equation}
holds with probability $q_1,$ where $0<q_1< q_0$.

We use the following fact to estimate the expected star discrepancy

\begin{equation}\label{EDNSX}
    \mathbb{E}[D^{*}_{N}(W)]=\int_{0}^{1}\mathbb{P}(D^{*}_{N}(W)\ge t)dt,
\end{equation}
where $D^{*}_{N}(W)$ denotes the star discrepancy of point set $W$.

Plugging $q_0$ into \eqref{dnsy1}, we have

\begin{equation}\label{DNsY}
D_{N}^{*}\left(Y\right)\leq (\sqrt{2}\cdot\Sigma_0+1)\frac{A(d,q_0,N)}{N}
\end{equation}
holds with probability $q_0$. Then \eqref{DNsY} is equivalent to

$$
    \mathbb{P}\big(D_{N}^{*}\left(Y\right)\ge (\sqrt{2}\cdot\Sigma_0+1)\frac{A(d,q_0,N)}{N}\big)=1-q_0.
$$

Now releasing $q_0$ and let

\begin{equation}\label{tsig0}
t=(\sqrt{2}\cdot\Sigma_0+1)\frac{A(d,q_0,N)}{N},
\end{equation}

\begin{equation}\label{c0sig0}
C_0(\Sigma_0,N)=\frac{\sqrt{2}\cdot\Sigma_0+1}{N},
\end{equation}
and
\begin{equation}\label{c1sig0}
C_1(d,\Sigma_0,N)=\frac{\sqrt{2}\cdot\Sigma_0+1}{N}\cdot(d\ln(2e)+d\ln(N+1)-\frac{\ln(2\pi d)}{2}).
\end{equation}

Then

\begin{equation}\label{tc1sig0}
    t=C_1(d,\Sigma_0,N)-C_0(\Sigma_0,N)\ln (1-q_0).
\end{equation}

Thus from \eqref{EDNSX} and $q_0\in[q,1)$, we have

\begin{equation}\label{ednstarz}
\begin{aligned}
&\mathbb{E}[D^{*}_{N}(Y)]=\int_{0}^{1}\mathbb{P}(D^{*}_{N}(Y)\ge t)dt\\&=\int_{1-e^{\frac{C_1(d,\Sigma_0,N)}{C_0(\Sigma_0,N)}}}^{1-e^{\frac{C_1(d,\Sigma_0,N)-1}{C_0(\Sigma_0,N)}}}\mathbb{P}\Big(D^{*}_{N}(Y)\ge (\sqrt{2}\cdot\Sigma_0+1)\frac{A(d,q_0,N)}{N}\Big)\cdot C_0(\Sigma_0,N)\cdot\frac{1}{1-q_0}dq_0%%\\&=\int_{q}^{1-e^{\frac{C_1(d,\Sigma_0,N)-1}{C_0(\Sigma_0,N)}}}\mathbb{P}\Big(D^{*}_{N}(Z)\ge (\sqrt{2}\cdot\Sigma_0+1)\frac{A(d,q_0,N)}{N}\Big)\cdot C_0(\Sigma_0,N)\cdot\frac{1}{1-q_0}dq_0
\\&=\int_{q}^{1-e^{\frac{C_1(d,\Sigma_0,N)-1}{C_0(\Sigma_0,N)}}}C_0(\Sigma_0,N)\cdot\frac{1-q_0}{1-q_0}dq_0.
\end{aligned}
\end{equation}

Furthermore, from \eqref{dnsxq1}, we have

$$
    \mathbb{P}\big(D_{N}^{*}\left(X\right)\ge (\sqrt{2}\cdot\Sigma_0+1)\frac{A(d,q_0,N)}{N}\big)=1-q_1.
$$

Following the steps from \eqref{tsig0} to \eqref{tc1sig0}, we obtain,

$$
    \mathbb{E}[D^{*}_{N}(X)]=\int_{0}^{1}\mathbb{P}(D^{*}_{N}(X)\ge t)dt=\int_{q}^{1-e^{\frac{C_1(d,\Sigma_0,N)-1}{C_0(\Sigma_0,N)}}}C_0(\Sigma_0,N)\cdot\frac{1-q_1}{1-q_0}dq_0.
$$

From $q_1< q_0,$ we obtain

$$
    \frac{1-q_1}{1-q_0}>\frac{1-q_0}{1-q_0}.
$$

Hence,

\begin{equation}\label{ewycom1}
    \mathbb{E}(D_{N}^{*}(Y))<\mathbb{E}(D_{N}^{*}(X)).
\end{equation}

\section{Conclusion}\label{conclu}

We study expected star discrepancy under jittered sampling. A strong partition principle for the star discrepancy is presented. The next most direct task is to find the expected star discrepancy better than jittered sampling, which may be closely related to the variance of the characteristic function defined on the convex test sets.

\end{document}